\documentclass[11pt,doublespace,vi]{article}
\usepackage{amsfonts}
\usepackage{subfigure} 
\usepackage{epsfig}
\usepackage{amsmath, amsthm}
\newcommand{\be}{\begin{eqnarray}}     	\newcommand{\ee}{\end{eqnarray}}

\newcommand{\vol}{\mathrm{Vol}}
\newcommand{\diam}{\mathrm{diam}}
\newcommand{\sing}{\mathrm{Sing}}
\newcommand{\ric}{\mathrm{Ric}}

\newcommand{\rem}{\mathrm{Rm}}
\newcommand{\dist}{\mathrm{dist}}
 
\title{Convergence of a K\"ahler-Ricci flow}
\author{Natasa Sesum}
\date{}
\theoremstyle{plain}
\newtheorem{dummy}{Dummy}

\theoremstyle{definition}

\newtheorem{remark}[dummy]{Remark}
\newtheorem{lemma}[dummy]{Lemma}
\newtheorem{theorem}[dummy]{Theorem}
\newtheorem{proposition}[dummy]{Proposition}
	
\newtheorem{definition}[dummy]{Definition}

\newtheorem{claim}[dummy]{Claim}

\begin{document}

\maketitle

\begin{abstract}
In this paper we prove that for a given K\"ahler-Ricci flow with
uniformly bounded Ricci curvatures in an arbitrary dimension, for
every sequence of times $t_i$ converging to infinity, there exists a
subsequence such that $(M,g(t_i + t))\to (Y,\bar{g}(t))$ and the
convergence is smooth outside a singular set (which is a set of
codimension at least $4$) to a solution of a flow. We also prove that
in the case of complex dimension $2$, without any curvature
assumptions we can find a subsequence of times such that we have a
convergence to a K\"ahler-Ricci soliton, away from finitely many
isolated singularities.
\end{abstract}

\begin{section}{Introduction}

Let $M$ be a compact K\"ahler manifold of dimension $n$ with the
K\"ahler metric $ds^2 = g_{i\bar{j}}dz^i d\bar{z}^j$. The Ricci
curvature of this metric is given by the formula

$$R_{i\bar{j}} = \frac{-\partial^2}{\partial
z^i\partial\bar{z}^j}\ln\det(g_{i\bar{j}}).$$
This implies that $\frac{\sqrt{-1}}{2\pi}R_{i\bar{j}}dz^i\and d\bar{z}^j$
is closed and its cohomology class is equal to the first 
Chern class $c_1(M)$ of $M$. We will assume that $c_1(M)$ is positive and that
it is represented by a K\"ahler form. We will consider the complex
version of Hamilton's Ricci flow equation of the following type

\begin{equation}
\label{equation-equation_KR_flow}
(g_{i\bar{j}})_t = g_{i\bar{j}} - R_{i\bar{j}} =
\partial_i\bar{\partial}_j u,
\end{equation}
where $g_{i\bar{j}}(t) = g_{i\bar{j}}(0) + \partial_i\bar{\partial}_j\phi$
and $\frac{d}{dt}\phi = u$.
In \cite{cap1985} H.D. Cao proved that a solution of
(\ref{equation-equation_KR_flow}) exists for all times $t\in
[0,\infty)$. A natural question that one can ask is what happens to a
flow when time approaches infinity. Under which conditions will it
converge?  How can we describe the objects that we get in a limit? In
this paper we will give partial answers to these questions.

In section $3$ we will consider a K\"ahler-Ricci flow
(\ref{equation-equation_KR_flow})  with uniformly bounded Ricci
curvatures. Our goal is to prove the following theorem.

\begin{theorem}
\label{theorem-theorem_KR_flow}
Assume we are given a flow (\ref{equation-equation_KR_flow}). Assume
that the Ricci curvatures are uniformly bounded, i.e. $|\ric| \le C$
for all $t$. Then for every sequence $t_i\to\infty$ there exists a
subsequence such that $(M,g(t_i + t))\to (Y,\bar{g}(t))$ and the
convergence is smooth outside a singular set $S$, which is at least of
codimension four. Moreover, $\bar{g}(t)$ solves the K\"ahler-Ricci
flow equation off the singular set.
\end{theorem}

In section $4$ we will restrict ourselves to $2$ dimensional complex
orbifolds, without any curvature assumptions. We want to prove the
following theorem.

\begin{theorem}
\label{theorem-theorem_soliton}
Let $g_{k\bar{j}}(t) = g_{k\bar{j}} - R_{k\bar{j}}$ be a
K\"ahler-Ricci flow on a $2$ dimensional complex, K\"ahler
manifold. Then for every sequence $t_i\to\infty$ there exists a
subsequence so that $(M,g(t_i+t))\to (Y,\bar{g}(t))$ and $\bar{g}(t)$
is a K\"ahler-Ricci soliton.
\end{theorem}

{\bf Acknowledgements:} The author would like to thank her advisor
Gang Tian for constant help, support and many useful discussions.

\end{section}

\begin{section}{Background and notation}

First of all, let us recall the definitions of Ricci solitons.

\begin{definition}
A solution $g_{i\bar{j}}$ to equation
(\ref{equation-equation_KR_flow}) on $M$ is called K\"ahler-Ricci
soliton if it moves along (\ref{equation-equation_KR_flow}) under
one-parameter family of automorphisms of $M$ generated by some
holomorphic vector field. 
\end{definition}
This means that 

$$g_{i\bar{j}} - R_{i\bar{j}} = V_{i,\bar{j}} + V_{\bar{j},i},$$ 
for some holomorphic vector field $V = (V^i)$. 
In the case of  limit solitons in Theorem
\ref{theorem-theorem_soliton}, we will show that the vector fields
come from the gradients of  functions on $M$, i.e. that

$$g_{i\bar{j}} - R_{i\bar{j}} = f_{,i\bar{j}},$$
and $f_{,ij} = 0$ for some real valued function $f$ on $M$.
This condition is equivalent to a fact that $V = \nabla f$
is a holomorphic vector field.

Perelman's functional $\mathcal{W}$ for a flow
(\ref{equation-equation_KR_flow}) is

$$\mathcal{W}(g,f,\tau) = (4\pi\tau)^{-n}\int_M e^{-f}
[2\tau(R + |\nabla f|^2) + f - 2n] dV_{g},$$
with a constraint that $(4\pi\tau)^{-n}\int_M e^{-f}dV_g =1$.
Perelman (\cite{perelman}) has proved some very interesting 
properties of flow (\ref{equation-equation_KR_flow}). We will
list them in the following theorem.

\begin{theorem}[Perelman]
\label{theorem-theorem_useful}
If (\ref{equation-equation_KR_flow}) is a flow on a complex, K\"ahler,
closed manifold $M$, then
\begin{enumerate}
\item
$C^{1,\alpha}$ norms of functions $u(t)$ are uniformly bounded along
the flow,
\item
the scalar curvatures $R(t)$ and the diameters $\diam(M,g(t))$ are
uniformly bounded along the flow,
\item
a volume noncollapsing condition holds along the flow, i.e.
there exists $C = C(g(0))$ such that $\vol_t(B(p,r)) \ge Cr^n$.
\end{enumerate}
\end{theorem}

We will need a theorem proved by Cheeger, Colding and Tian
(\cite{tian2002}) in our further discussion and we will state it below
for a reader's convenience.

\begin{theorem}[Cheeger, Colding, Tian]
\label{theorem-theorem_cct}
If $\{M_i,g_i,p_i)\}$ converges to $(Y,d,y)$ in pointed
Gromov-Hausdorff topology, if $|\ric|_{M_i} \le C$ and if
$\vol(B_1(p_i)) \ge C$ for all $i$, then the regular part
$\mathcal{R}$ of $Y$ is a $C^{1,\alpha}$-Riemannian manifold and at
points of $\mathcal{R}$, the convergence is $C^{1,\alpha}$. Moreover
the codimension of the set of singular points (which is a closed set
in $Y$) is at least $4$.
\end{theorem}

In the proof of Theorem \ref{theorem-theorem_KR_flow} we will use
Perelman's pseudolocality theorem (\cite{perelman2002}).

\begin{theorem}[Perelman]
\label{theorem-theorem_pseudolocality}
For every $\alpha > 0$ there exist $\delta > 0$, $\epsilon > 0$ with
the following property. Suppose we have a smooth solution to the Ricci
flow and assume that at $t=0$ we have $R(x) \ge -r_0^{-2}$ and
$\vol(\partial\Omega)^n \ge (1-\delta)c_n\vol(\Omega)^{n-1}$ for any
$x,\Omega \subset B(x_0,r_0)$, where $c_n$ is the euclidean
isoperimetric constant. Then, $|\rem|(x,t) \le \alpha t^{-1} +
(\epsilon r_0)^{-2}$ whenever $0 < t \le (\epsilon r_0)^2$ and
$\dist_t(x,x_0) < \epsilon r_0$.
\end{theorem}
Perelman proved this theorem for a case of unnormalized Ricci flow,
but it can be easily modified for the case of a normalized
K\"ahler-Ricci flow.

\end{section}

\begin{section}{K\"ahler-Ricci flow with uniformly bounded Ricci curvatures}

In this section we will consider a flow
(\ref{equation-equation_KR_flow}), with uniformly bounded Ricci
curvatures. For any sequence $t_i\to\infty$, if $g_i(t) = g(t_i+t)$,
the metrics $g_i(t)$ are uniformly equivalent to metrics $g_i(s)$ for
$s,t$ belonging to an interval of finite length. Moreover, the
following proposition (in \cite{glickenstein2002}) applies to metrics
$g_i(t)$.

\begin{proposition}[D. Glickenstein]
\label{proposition-proposition_parameter_convergence}
Let $\{(M_i,g_i(t),p_i)\}_{i=1}^{\infty}$, where $t\in[0,T]$, be a
sequence of pointed Riemannian manifolds of dimension $n$ which is
continuous in the $t$ variable in the following way: for each $\delta
> 0$ there exists $\eta > 0$ such that if $t_0,t_1\in [0,T]$ satisfies
$|t_0 - t_1| < \eta$ then

\begin{equation}
\label{equation-equation_condition}
(1 + \delta)^{-1}g_i(t_0) \le g_i(t_1) \le (1 + \delta)g_i(t_0),
\end{equation}
for all $i > 0$, and such that $\ric(g_i(t)) \ge cg_i(t)$, where $c$
does not depend on $t$ or $i$. Then there is a subsequence
$\{(M_i,g_i(t),p_i)\}_{i=1}^{\infty}$ and a $1$-parameter family of
complete pointed metric spaces $(X(t),d(t),x)$ such that for each
$t\in [0,T]$ the subsequence converges to $(X(t),d(t),x)$ in the
pointed Gromov-Hausdorff topology.
\end{proposition}

$(M_i,g_i(t))$ and $(M_i,g_i(0))$ are homeomorphic by Lipschitz
homeomorphisms, and in \cite{glickenstein2002} it has been showed be
showed that $X(t)$ is homeomorphic to $X(0)$. If $t_i$ is any sequence
such that $t_i\to\infty$, Proposition
\ref{proposition-proposition_parameter_convergence}  applies to
$(M,g(t_i + t))$ for all $i$ and all $t$ belonging to a time interval
of finite length.

For the moment we will restrict ourselves to the case of K\"ahler
manifolds of complex dimension $2$, and later we will show how it can
be generalized to an arbitrary dimension. In the case of complex
dimension $2$, for every sequence $t_i\to\infty$ there is a
subsequence $\{(M,g(t_i + t))\}$ converging to a compact orbifold
$(Y,\bar{g}(t))$ with isolated singularities. This is due to the fact
that $L^2$ norm of the curvature operator in the K\"ahler case can be
uniformly bounded in terms of the first and the second Chern class of
a manifold and its K\"ahler class. Combining Proposition
\ref{proposition-proposition_parameter_convergence} and Theorem
\ref{theorem-theorem_cct} gives that $(Y,\bar{g}(t))$ is an $1$
parameter family of orbifolds (it is even a Lipschitz family for $t$
belonging to an interval of finite length), such that a regular part
of $(Y,\bar{g}(t))$ is $C^{1,\alpha}$ manifold and the convergence
$(M,g(t_i + t))\to (Y,\bar{g}(t))$ takes place in $C^{1,\alpha}$
topology, away from the set of singular points (which is common for
all orbifolds $(Y,\bar{g}(t))$). In the case of higher dimensions,
again by Proposition
\ref{proposition-proposition_parameter_convergence} and Theorem
\ref{theorem-theorem_cct} we will have that $\{(M,g(t_i+t))\}$
converge to $(Y,\bar{g}(t))$ with a singular set $S\subset Y$ of
codimension at least $4$. $\mathcal{R} = Y\backslash S$ is an open
$C^{1,\alpha}$ manifold and the convergence on $\mathcal{R}$ is in
$C^{1,\alpha}$ norm. We will show later that the set $S$ is common for
$(Y,\bar{g}(t))$ for all $t$. The main tools in the proof of Theorem
\ref{theorem-theorem_KR_flow} will be Theorem
\ref{theorem-theorem_pseudolocality} and Theorem $A.1.5$ of Cheeger
and Colding  that can be found in the appendix of
\cite{cheeger1997}.

We will now prove Theorem \ref{theorem-theorem_KR_flow}.

\begin{proof}

If the curvature does not blow up, we are done. Therefore, assume that
the curvature does blow up. Let $t_i\to\infty$ be such that $Q_i =
|\rem|(p_i,t_i) \ge \max_{M\times [0,t_i]}|\rem|(x,t)$ and
$Q_i\to\infty$.  We already know that since $|\ric|(g(t)) \le C$,
there exists a subsequence $(M,g(t_i + t))$ converging to orbifolds
$(Y,\bar{g}(t))$ in $C^{1,\alpha}$ norm off the set of singular
points. Moreover, metrics $\bar{g}(t)$ are $C^{1,\alpha}$ off the
singular set. We may assume that $\sing(Y) = \{p\}$. Our goal is to
show that we actually have $C^{\infty}$ convergence off the singular
point $p$, due to the fact that our metrics are changing with the
K\"ahler-Ricci flow.

Adopt the notation of \cite{tian2002}. In general, a point $y\in Y$ is
called regular, if for some $k$, every tangent cone at $y$ is
isometric to $\mathrm{R}^k$. Denote a set of those points by
$\mathcal{R}_k$ and let $\mathcal{R} = \cup_k\mathcal{R}_k$. Because
of the noncollapsing condition that we have because of Theorem
\ref{theorem-theorem_useful}, we have that $\mathcal{R} =
\mathcal{R}_n$. Let $\mathcal{R}_{\epsilon} = \{y\:\:|\:\:
d_{GH}(B_1(y_{\infty}), B_1(0)) < \epsilon$ for every tangent cone
$(Y_y,y_{\infty}) \}$ , where $B_1(0)$ is a unit ball in
$\mathrm{R}^n$. Let $\mathcal{R}_{\epsilon,r}$ be a set of all
points $y\in Y$ such that there exists $x$ such that $(0,x)\in
R^4\times \{x\}$ and for some $u > r$ and every $s\in (0,u]$
$d_{GH}(B_s(y),B_s((0,x))) < \epsilon s$. $\mathcal{R}_{\epsilon} =
\cup_r \mathcal{R}_{\epsilon, r}$.

Choose $\epsilon_P$ and $\delta_P$ as in Perelman's pseudolocality
theorem. Choose $\epsilon' > 0$ such that $\delta_P > \epsilon'$ and
$\epsilon' \le \epsilon_0$, where $\epsilon_0$ is such that
$\mathcal{R} = \mathcal{R}_{\epsilon}$ for all $\epsilon \le
\epsilon_0$ (the existence of such an $\epsilon_0$ is proved in
section $7$ of \cite{cheeger1997}.

Fix a time $t = t_0$. Pick up any point $q\in Y\backslash\{p\}$. Then
$q\in \cap_{\epsilon \le \epsilon_0} \mathcal{R}_{\epsilon}$. Let $d =
\dist_{\bar{g}(t_0)}(p,q)$.

\begin{claim}
\label{claim-claim_convergence}
There exist $\eta > 0$ and a sequence $q_i \in M$ such that $q_i\to
q$, while $(M,g(t_i+t))\to(Y,\bar{g}(t))$ as $i\to\infty$ for all $t$.
\end{claim}

\begin{proof}
Assume that $q\in K \subset Y\backslash\{p\}$, where $K$ is a compact
set and let $r = \dist (K,p)$. For every $t$, $g(t_i+t)$ uniformly
converge to $\bar{g}(t)$ on $K$. Let $\phi: K\to K_i$ be
diffeomorphisms as in a definition of convergence of $(M,g(t_i+t))$ to
$(Y,\bar{g}(t))$. Let $q_i(t_0)\in M$ be such that
$\dist_{g(t_i+t_0)}(q_i(t_0),\phi_i(q)) < \epsilon$, for $i\ge
i_0$. Since the Ricci curvatures of $g(t)$ are uniformly bounded, there
exists $\eta > 0$ so that $|t - s| < \eta$ implies that
$|\dist_{g(t_i+s)}(x,y) - \dist_{g(t_i+t)}(x,y)| < \epsilon$, for all
$x,y \in M$. Therefore,

$$\dist_{g(t_i+t)}(q_i(t_0),\phi_i(q)) \le
\dist_{g(t_i+t_0)}(q_i(t_0),\phi_i(q)) + \epsilon < 2\epsilon,$$ 
for $i\ge i_0$ and for all $t\in [-\eta,\eta]$.  Notice that $\eta$ does
not depend either on $t_0$ or $q$, but it depends on $K$, i.e.  on its
distance from $p$. Therefore, if we continue this process infinitely
many times, considering $t_0+\eta$ instead of $t_0$, etc. we  get that 
the sequence $\{q_i\}$ will work for all times $t\ge 0$.

\end{proof}

\begin{lemma}
\label{lemma-lemma_perelman}
For any regular point $q\in \mathcal{R}$ there exists $i_0$, $\eta$
and $r > 0$ such that for all $B_{g(t_i+t)}(s,q')\subset
B_{g(t_i+t)}(r,q_i(t))$ we have $\vol_{g(t_i+t)}B_{g(t_i+t)}(s,q') \ge
(1 - \epsilon')s^n$, for all $i \ge i_0$ and all $t\in
[t_0-\eta,t_0+\eta]$, where $q_i\in M$ is a sequence of points
converging to $q$, while $(M,g(t_i + t))\to (Y,\bar{g}(t))$.
\end{lemma}

\begin{proof}

For $\epsilon'$ find $r$ and $\delta$ as in Theorem $A.1.5$ (i) and
(ii) in \cite{cheeger1997}. For this $\delta$ (that now plays the role
of $\epsilon$ in Theorem $A.1.5$ in \cite{cheeger1997}) find
$\delta_1$ and $r_1$ (by part (iii) of the same theorem), such that
$x\in (\mathcal{W}\mathcal{R})_{8\delta_1,r'}$ implies that

\begin{equation}
\label{equation-equation_e_delta}
y\in \mathcal{R}_{\delta,s}  \:\:\: \forall y\in B_{r'}(x) \:\:\:\forall
s \le (1 - \delta)r' - \dist_0(x,y),\:\:\: r' \le r_1,
\end{equation}
where a distance is measured in metric $\bar{g}(t_0)$. We may assume
that $r_1 < d$, because otherwise we can decrease $r_1$. Take any
sequence $\delta_i\to 0$ as $i\to\infty$. We can choose a sequence
$r_i$ such that $q\in \mathcal{R}_{\delta_i,r_i}$, since $q\in
\mathcal{R}_{\delta_i}$. We claim that $q\in \mathcal{R}_{\delta_1,r}$,
for some $r < r_1$. In order to prove that, we may assume $r_i\to 0$
(otherwise if $r_i \ge \kappa$ for all $i$, $d_{GH}(B_l(q), B_l(0))
\le l\delta_i \to 0$ for all $l \le \kappa$ and therefore we would have
$\vol B_{l}(q) = \vol B_{l}(0)$ for all $l \le \kappa$, and $q\in
\mathcal{R}_{\delta_1,s}$ for some $s < r_1$, by Theorem $A.1.5$, part
(i) in \cite{cheeger1997}. Therefore, there exist $\delta'' <
\delta_1$ and $r'' < r_1$ such that $q\in
\mathcal{R}_{\delta'',r''}$. This implies $q\in
\mathcal{R}_{\delta_1,r''}$, since $\delta'' < \delta_1$.  This is
true in metric $\bar{g}(t_0)$.

\begin{claim}
\label{claim-claim_small}
There exist $\eta > 0$ and $i_0$ such that $q_i\in
\mathcal{R}_{\delta_1,r''}$ for all metrics $g(t_i+t)$ for $i\ge i_0$
and $t\in [t_0-\eta,t_+\eta]$.
\end{claim}

\begin{proof}

$q\in \mathcal{R}_{\delta_1,r''}$ and therefore, 

\begin{equation}
\label{equation-equation_eq3}
d_{GH}(B_s(x,0), B_s(q,t_0)) < s\delta_1,
\end{equation}
for some $s < r''$. We
can substitute space $\{x\}\times \mathrm{R}^4$ by $\mathrm{R}^4$ only
and therefore we can write just $B_s(0)$ instead of $B_s(x,0)$.  Since
the Ricci curvatures of $g(t)$ are uniformly bounded, there exists
$\eta$ such that $|t-t_0| < \eta$ implies that
\begin{equation}
\label{equation-equation_eq1}
d_{GH}(B_{g(t_i+t)}(q_i,s),B_{g(t_i+t_0)}(q_i,s)) < \delta_1 s.
\end{equation}
Since $g(t_i+t_0)$ converges to $\bar{g}(t_0)$ uniformly, away from a
singular point $p$, there exists $i_0$ (depending on $\delta_1 s$ and
a compact set $K$) such that for $i\ge i_0$
\begin{equation}
\label{equation-equation_eq2}
d_{GH}(B_{g(t_i+t_0)}(q_i,s), B_{\bar{g}(t_0)}(q,s)) < \delta_1 s.
\end{equation}

Combining estimates (\ref{equation-equation_eq3}),
(\ref{equation-equation_eq1}) and (\ref{equation-equation_eq2}),
together with an approximate triangle inequality for Gromov-Hausdorff 
distance we get

$$d_{GH}(B_{g(t_i+t)}(q_i,s), B_s(0)) < 4\delta_1 s,$$
for all $i\ge i_0$ and all $t\in[t_0-\eta,t_0+\eta]$. This implies that
$q_i\in\mathcal{W}\mathcal{R}_{8\delta_1,r''}$, for all $i\ge i_0$
and all $t\in [t_0-\eta,t_0+\eta]$.

\end{proof}

Combining Claim \ref{claim-claim_small} and part (iii) of Theorem
$A.1.5$ in \cite{cheeger1997}, we get that $q'\in
\mathcal{R}_{\delta,s}$, for all $q'\in B_{g(t_i+t)}(q_i,r'')$, $s\le
(1-\delta)r'' - \dist_{t_i+t}(q_i,q')$, for all $i\ge i_0$ and $t\in
[t_0-\eta,t_0+\eta]$. Part (ii) of Theorem $A.1.5$ in
\cite{cheeger1997} gives that

\begin{equation}
\label{equation-equation_volume_inner}
\vol B_{g(t_i+t)}(s,q') \ge (1 - \epsilon')\vol B_s(0),
\end{equation}
for all $q'\in B_{g(t_i+t)}(r'',q)$ and $s \le (1 - \delta)r'' -
\dist_{g(t_i+t)}(q_i,q')$. By reducing $r''$ we get that there 
exists $r''$ such that the estimate
(\ref{equation-equation_volume_inner}) holds for all $q'\in
B_{g(t_i+t)}(r'',q)$ and all $s$ such that $B_{g(t_i+t)}(s,q')\subset
B_{g(t_i+t)}(r'',q)$, for $i\ge i_0$ and $t\in [t_0-\eta,t_0+\eta]$.

\end{proof}

Choose $r$, $i_0$ and $\eta$ as in the claim above (for our regular
point $q$ that we have fixed earlier). Reduce $r''$ if necessary, so
that $(\epsilon'r'')^2 < \eta$. Since $1 - \epsilon' > 1 - \delta_P$,
and since for every ball $B_{g(t_i-(\epsilon'r'')^2/2)}(q',s)\subset
B_{g(t_i-(\epsilon'r'')^2/2)}(q_i,r'')$, we have that
$\vol_{g(t_i-(\epsilon'r'')^2/2)}B_s(q') \ge (1 - \delta_P)s^nc_n$, by
Perelman's pseudolocality Theorem \ref{theorem-theorem_pseudolocality}

$$|\rem|(x,t) \le \frac{1}{(\epsilon'r'')^2} + (\epsilon'r'')^2,$$ for
all $x\in B_{g(t)}(q_i,\epsilon'r'')$ and for every $t\in [t_i -
(\epsilon'r'')^2/2, t_i + (\epsilon'r'')^2/2]$. We have that $g_i(t) =
g(t_i + t)$ is a sequence of Ricci flows with uniformly bounded
curvatures for $t\in [-(\epsilon'r'')^2/2, (\epsilon'r'')^2/2]$ on
balls $B_{g_i(t)}(q_i,\epsilon'r'')$. This together with the volume
noncollapsing condition and Hamilton's compactness theorem give that
the convergence of the sequence of our metrics is smooth, and
$\bar{g}(t)$ are smooth metrics on $B_{\bar{g}(t)}(q,\epsilon'r'')$,
for $t\in [0,(\epsilon'r'')^2/2]$. Repeating the procedure described
above infinitely many times, to time intervals translated by
$(\epsilon'r'')^2/2$ (considering $t_0+(\epsilon'r'')^2/2$ instead of
$t_0$, etc.) and applying diagonalization method to a sequence of
times $t_i$ (since for every step of length $(\epsilon'r'')^2/2$ we
have to extract a subsequence of a subsequence), we get that
$\bar{g}(t)$ are smooth metrics on $B(q,\epsilon'r'')$ for all times
$t \ge 0$ (we can take $t_0 = 0$) and that $g(t_i + t) \to \bar{g}(t)$
smoothly on $B_{\bar{g}(t)}(q,\frac{\epsilon'r''}{2})$ for all times
$t \ge 0$. We will use the fact that the Ricci tensor is uniformly
bounded to show that we can extend the previous result from a ball to
any compact set $K\subset Y$. By a definition of convergence, that
will mean $(M,g(t_i + t)) \to (Y,\bar{g}(t))$ smoothly, away from the
set of singular points.

Take a compact set $K\subset Y\backslash S$, where $S$ is a set of
singular points on $(Y,\bar{g}(t))$. It is the same set for all
singular metrics $\bar{g}(t)$. Let $\phi_i: K\to K_i$ be a sequence of
diffeomorphisms from a definition of convergence of metrics $g_i(0)$
to a metric $\bar{g}(0)$.

$|\ric|(t) \le C$ for all $t$ by the assumtion of the main theorem. We
have proved that $\bar{g}(t)$ is 1-parameter family of metrics on
$Y$. Moreover, $\bar{g}(t)$ satisfies the K\"ahler-Ricci flow equation
away from the singular points.

\begin{claim}
There exist $\delta > 0$, a subsequence $t_i$ and $C_1 = C_1(K)$ such
that $|\rem|(g(t_i+t)) \le C$ on $K_i$ for all $t\in [t_0,t_0\delta]$
and all $t_0$.
\end{claim}

\begin{proof}
Fix $t_0$.  For every $q\in Y$ we can choose $r_q > 0$, $\eta_q$ and
$i_q$ as in Lemma \ref{lemma-lemma_perelman}. Look at the collection
of balls $B_{\bar{g}(t_0)}(q,(\epsilon'r_q)/4)$ covering $K$. Since
$K$ is compact we can consider only finitely many of them covering
$K$. Denote their centres by $q_1,q_2,\dots q_N$.  Since $\bar{g}(t)$
solves the equation (\ref{equation-equation_KR_flow}) and since the Ricci
curvatures of $\bar{g}(t)$ are uniformly bounded on
$Y\backslash\{p\}$, there exists $A>0$ so that the balls
$B_{\bar{g}(t)}(q_i,(\epsilon'r_{q_i})/2)$ cover $K$, for $t\in
[t_0-A,t_0+A]$. Let $r_1 = \min\{r_{q_1},r_{q_2},\dots,r_{q_N}\}$,
$\eta_1 = \min\{{\eta}_{q_1}, \dots, {\eta}_{q_N}\}$ and $i_1 =
\max\{i_{q_1},\dots i_{q_N}\}$. Then $|\rem| (x,t + t_i) \le
\frac{1}{(\epsilon'r_1)^2} + (\epsilon'r_1)^2 = C_1(K)$ for all $x\in
B_{g(t_i+t)}(q_i^j,\epsilon' r_{q_j})$, all $i\ge i_1$ and all $t\in
[0,\min\{{\eta}_1,(\epsilon'r_1)^2/2\}]$, where $q_i^j$ are the
sequences of points such that $B_{g_(t_i+t)}(q_i^j,r_{q_j}) \to
B_{\bar{g}(t)}(q_j,r_{q_j})$ while
$(M,g(t_i+t))\to(Y,\bar{g}(t))$. Let $\delta =
\min\{{\eta}_1,(\epsilon'r_1)^2/2,A\}]$.

The balls $B_{\bar{g}(t)}(q_j,(\epsilon'r_{q_j})/2)$ for $1\le j\le N$
cover $K$. A definition of convergence gives that there exists $i_0
\ge i_1$ so that $B_{g(t_i)}(q_i^j,\frac{2\epsilon'r_{q_j}}{3})$ cover
$K_i$ for all $i \ge i_0$. We can assume that $\delta$ is small enough
so that

$$B_{g(t_i)}(q_i^j,(r_{q_j}-a)\epsilon') \subset
B_{g(t_i+s)}(q_i^j,r_{q_j}\epsilon'/2),$$ for $a < \frac{r_1}{3}$ so
that $r_{q_j} - a > r_{q_j} - \frac{r_1}{3} > \frac{2r_{q_j}}{3}$ and
therefore $B_{g(t_i)}(q_i^j,(\epsilon'(2r_{q_j})/3)\subset
B_{g(t_i)}(q_i^j, (r_{q_i^j} - a)\epsilon')$. Since the balls
$B_{g(t_i)}(q_i^j,(2r_{q_j}\epsilon')/3)$ cover $K_i$, so do balls
$B_{g(t_i+s)}(q_i^j,r_{q_j})$ for all $i\ge i_0$ and all
$s\in[t_0-\delta,t_0+\delta]$. Therefore, $|\rem|(x,t_i+s) \le C_1$
for all $i\ge i_1$ and all $x\in K_i$ and all $s\in
[t_0-\delta,t_0+\delta]$. Therefore we actually can extract a
subsequence $t_i$ such that the pullbacks of metrics $g(t_i+s)$
converge to a solution of the K\"ahler-Ricci flow uniformly on
$K\times [t_0-\delta,t_0+\delta]$.
\end{proof}

Applying the method from the previous claim to a sequence $t_i+\delta$
instead of a sequence $t_i$ we can find a subsequence such that
$g(t_i+s)\to\bar{g}(t)$ smoothly and uniformly on
$K\times[t_0,t_0+2\delta]$, since our choice of $\delta$ does not
depend on an initial time, but on a chosen compact set $K\subset Y$
and a uniform bound on the Ricci tensor. Repeating this infinitely
many times and diagonalizing the sequence $t_i$, we get a sequence
$t_i$ such that $g(t_i+s)\to\bar{g}(t)$ smoothly converge on all
compact subsets of $K\times [0,\infty)$. We can choose a countable
sequence of compact sets $L_k$ exhausting $Y\backslash S$. We can find
a subsequence of $t_i$ for each $L_k$, so that the above that we have
proved for any compact set $K\subset Y\backslash S$, applies to $L_k$
as well. By a diagonalization procedure applied to $t_i$ we can get a
subsequence so that $(M,g(t_i+t)\to (Y,\bar{g}(t))$ for all $t \ge 0$,
where $\bar{g}(t)$ is a solution to the K\"ahler-Ricci flow away from
the set of singular points. The convergence is in the sense that for
every compact set $K\subset Y\backslash S$ there exist diffeomorphisms
$\phi_i: K\to K_i$, where $K_i\subset M$ are compact and $\phi_i^*
g(t_i+t)\to\bar{g}(t)$, uniformly and smoothly on all compact subsets
of $K\times [0,\infty)$.

\end{proof}

In the proof of Theorem \ref{theorem-theorem_KR_flow} we assumed that
the complex dimension was $2$. The proof above generalizes to an
arbitrary dimension easily. Since the Ricci tensor is uniformly
bounded along the flow, for every sequence $t_i\to\infty$ there exists
a subsequence so that $(M,g(t_i + t)) \to (Y(t),\bar{g}(t))$ and the
convergence is smooth outside a set $S(t)$ of codimension at least
$4$. As above, it easily follows that $Y = Y(t)$ for all $t$. We
should only check that $S(t) = S(s)$ for any $s,t\in [0,\infty)$.

\begin{lemma}
\label{lemma-lemma_adjustion}
$S(s) = S(t)$ for any $s,t\in [0,\infty)$.
\end{lemma}

\begin{proof}

It is enough to prove: $\exists a > 0$ such that for $|s - t| < a$
$S(s) = S(t)$.

Choose $\epsilon > 0$ such that $\mathcal{R}_{\epsilon}(s) =
\mathcal{R}(s)$ and $\mathcal{R}_{\epsilon}(t) = \mathcal{R}(t)$, for
$|s - t| < a$, where we will choose $a$ later. Assume there exists
$q\in S(t)\backslash S(s)$. That implies $q\in \mathcal{R}(s)$.  For
$\epsilon > 0$ choose $\epsilon' = \epsilon'(\epsilon,n) > 0$ and $r'
= r'(\epsilon, n)$ so that Theorem $A.1.5$ in \cite{cheeger1997}
holds. Then the following claim holds for $q$.

\begin{claim}
There exist $i_0$ and $r < r'$ such that for all
$B_{g(t_i+s)}(q',u)\subset B_{g(t_i+s)}(q_i,r)$ we have
$\vol_{g(t_i+s)}B_{g(t_i+s)}(q',u) \ge c_n(1 - \epsilon'/2)u^n$, for
all $i \ge i_0$, where $q_i\in M$ is a sequence of points converging
to $q$, while $(M,g(t_i + s))\to (Y,\bar{g}(s))$.
\end{claim}

The proof of this claim is the same as the proof of Lemma
\ref{lemma-lemma_perelman}.

Since the Ricci tensors are uniformly bounded along the flow, we have
a good control on the volumes and the sizes of balls in metrics at
different times, when the considered time interval is sufficiently
small. Similarly as in \cite{natasa2003} we can find sufficiently
small $a > 0$ such that $|s - t| < a$, for any $u < r$ implies that

$$\vol_{g(t_i+t)}B_{g(t_i+t)}(q,u) \ge
\sqrt{(1-\frac{\epsilon'}{2})}\vol_{g(t_i+s)}B_{g(t_i+s)}(q,u) \ge (1
- \frac{\epsilon'}{2}),$$
$$B_{g(t_i+s)}(q,u\tilde{r}) \subset B_{g(t_i+t)}(q,u),$$ where
$\tilde{r} = \frac{1}{1+(e^{2C|s-t|}-1)^{\frac{1}{2}})}$ and we can
choose $a$ small enough, so that $\tilde{r}^n >
\sqrt{1-\frac{\epsilon'}{2}}$. Finally, since
$\vol_{g(t_i+s)}B_{g(t_i+s)}(q,u) \ge c_n(1-\frac{\epsilon'}{2})u^n$,
we get that

$$\vol_{g(t_i+t)}B_{g(t_i+t)}(q,u) \ge (1-\epsilon'/2)^2u^nc_n \ge
(1-\epsilon')c_nu^n,$$ i.e. $q\in \mathcal{R}_{\epsilon,r/2} \subset
\mathcal{R}_{\epsilon} = \mathcal{R}(t)$. This means that $q$ can not
be in $S(t)$ and we get a contradiction.  We can repeat the procedure
above infinitely many times to get that $S(t) = S(s)$ for all $s,t\in
[t_0,t_0+a]$ and all $t_0 \ge 0$, i.e. $S(t) = S$ for all $t\ge 0$.

Having Lemma \ref{lemma-lemma_adjustion} we can repeat the proof of
Theorem \ref{theorem-theorem_KR_flow} for complex dimension $2$, to
get that theorem is actually true for all dimensions.

\end{proof}

\end{section}

\begin{section}{K\"ahler-Ricci soliton as a limit}

Assume that $(Y,\bar{g})$ is a complex $n$-dimensional orbifold with
finitely many singularities. Assume without loss of generality that
$p$ is its only singular point. Analogously to the case of compact
manifolds we can show that $\mu(\bar{g}(t),\frac{1}{2}) =
\inf_{\{f|(2\pi)^{-n}\int_M e^{-f} =1\}}$ is achieved and that

$$2\Delta f - |\nabla f|^2 + R + f - 2n =
\mu(\bar{g}(t),\frac{1}{2}),$$ on $Y\backslash \{p\}$ where $f(t)$ is
a function such that $\mu(\bar{g}(t),\frac{1}{2}) =
W(\bar{g}(t),f(t),\frac{1}{2})$. Before we start proving Theorem
\ref{theorem-theorem_soliton}, we will fisrt prove the following
proposition.

\begin{proposition}
\label{proposition-proposition_mu}
Let $(M,g_i)$ be a sequence of smooth, closed manifolds, with
uniformly bounded scalar curvatures, converging to an orbifold
$(Y,\bar{g})$ with a singular point $p$. Assume that
$|\mu(g_i,\frac{1}{2})| \le C$ for all $i$. Then
$\lim_{i\to\infty}\mu(g_i,\frac{1}{2}) = \mu(\bar{g},\frac{1}{2})$.
\end{proposition}

\begin{proof}

Fix $\epsilon > 0$. Similarly as in the smooth case we can show that
$|f_i|_{C^{2,\alpha}} \le C$ on $Y\backslash\{p\}$ (by weak regularity
theory applied to $\tilde{f}_i = e^{-\frac{f_i}{2}}$ we get that
$W^{3,p}$ norms of $\tilde{f_i}$ are uniformly bounded and Sobolev
embedding theorems applied to $\tilde{f_i}$ and $\Omega_k = Y\backslash
B(p,\frac{1}{k})$ give uniform $C^{2,\alpha}$ bounds,
i.e. $|\tilde{f_i}|_{C^{2,\alpha}(\Omega_k)} \le C$, where $C$ does not
depend on $k$).

$\mu(\bar{g},\frac{1}{2}) = \int_Y F d_{\bar{g}}$, where $F =
(2\pi)^{-n}e^{-f}(|\nabla f|^2 + R + f - 2n)$. $F$ is an integrable
function and therefore there exists some $r > 0$ such that
$\int_{B(p,2r)} |F| dV_{\bar{g}} < \epsilon$. Then:

\begin{eqnarray*}
\mu(\bar{g},\frac{1}{2}) &>& \int_{Y\backslash B(p,2r)} F dV_{\bar{g}} 
- \epsilon\\
&=& \int_{Y\backslash B(p,2r)} F(dV_{\bar{g}} - dV_{\tilde{g}_i}) +
\int_{Y\backslash B(p,2r)} F dV_{\tilde{g}_i} - \epsilon,
\end{eqnarray*}
for $\tilde{g}_i = \phi_i^* g_i$ ($\phi_i:Y\backslash B(p,2r) \to
M\backslash B(p_i, r)$ are diffeomorphisms from a definition of
convergence). Since $|F| \le C$ on $Y\backslash \{p\}$ and since
$\tilde{g}_i$ converge uniformly to $\bar{g}$ on $Y\backslash
B(p,2r)$, there exists $i_0$ such that for all $i \ge i_0$

\begin{eqnarray*}
\mu(\bar{g},\frac{1}{2}) &>& \int_{Y\backslash B(p,2r)} F
dV_{\tilde{g}_i} - 2\epsilon\\
&=& \int_{\phi_i(Y\backslash B(p,2r))} \phi^*F dV_{g_i} - 2\epsilon \\
&=& \int_{M\backslash U_r^i} \phi^*_i F dV_{g_i} - 2\epsilon,
\end{eqnarray*}
where $p_i\in U_r^i$ and $B(p_i,r)\subset U_r^i \subset B(p_i,2r)$ (we
can always choose big $i_0$ such that these inclusions hold for $i \ge
i_0$). Let $\tilde{f}_i = \phi_i^* f$. $\int_Y e^{-f} dV_{\bar{g}} =
(2\pi)^n$ and therefore $(2\pi)^n - \epsilon < \int_{Y\backslash
B(p,2r)} e^{-f} dV_{\bar{g}} < (2\pi)^n + \epsilon$ (we can decrease
$r > 0$ so that this holds, since $f$ is bounded on
$Y\backslash\{p\}$). For the same reasons as above, for $i\ge i_0$ we
have that $(2\pi)^{n} - 2\epsilon < \int_{Y\backslash B(p,2r)} e^{-f}
dV_{\tilde{g}_i} < (2\pi)^n + \epsilon$. Therefore, for $i \ge i_0$
(we increase $i_0$ if necessary)

$$(2\pi)^n - \epsilon < \int_{\phi_i(Y\backslash
B(p,2r))}e^{-\tilde{f}_i} dV_{g_i} < (2\pi)^n + 2\epsilon.$$ Choose
$\eta_i = 1$, outside $B_{g_i}(p_i,3r)$ and $0$ inside the ball
$B_{g_i}(p_i,2r)$ so that $|\nabla_i \eta_i| \le \frac{1}{r}$. Let
$\bar{f}_i = \eta_i \tilde{f}_i$. For $r$ small enough and $i_0$ big
enough we have

$$|\int_{\phi_i(Y\backslash B(p,2r))} e^{-\tilde{f}_i} dV_{g_i} -
\int_{\phi_i(Y\backslash B(p,2r))} e^{-\bar{f}_i} dV_{g_i}| <
\epsilon,$$ which implies $(2\pi)^n - C\epsilon < \int_M
e^{-\bar{f}_i} dV_{g_i} < (2\pi)^n + C\epsilon$. Modify each
$\bar{f}_i$ by a small constant so that $\int_M e^{-\bar{f}_i}
V_{g_i} = 1$. Furthermore,

\begin{eqnarray*}
& & |\int_M |\nabla\bar{f}_i|^2 e^{-\bar{f}_i} dV_{g_i} - \int_M
|\nabla\tilde{f}_i|^2 e^{-\tilde{f}_i} dV_{g_i}| \le \\
&\le& C\int_{B(p_i,3r)\backslash B(p_i,2r)} (|\nabla\eta_i|^2 \tilde{f}_i^2
+ \eta_i^2|\nabla \tilde{f}_i|^2 + |\nabla\tilde{f}_i|^2) dV_{g_i}\\ 
&\le& C(\int_{B(p_i,3r)\backslash B(p_i,2r)}\frac{1}{r^2}dV_{g_i} + 
\int_{B(p_i,3r)\backslash B(p_i,2r)}|\nabla\tilde{f}_i|^2 dV_{g_i}) 
< \epsilon,
\end{eqnarray*}
for small enough $r$, since there exists constant $C$ so that
$\vol_{g(t)}B(x,r) \le Cr^4$ for all $r \le r_0$ and all $t$. This
simple follows from the fact that $\frac{\vol B(p,r)}{V_{-C}(r)} \le
\frac{\vol B(p,\delta)}{V_{-C}(\delta)}$, where $V_{-C}(r)$ is a
volume of a ball of radius $r$ in a simply connected space of constant
curvature $-C$, the fact that $\lim_{\delta\to 0}\frac{\vol
B(p,\delta)}{V_{-C}(\delta)} = w_n$ and $\lim_{r\to 0}
\frac{V_{-C}(r)}{r^4} = w_n$ ($w_n$ is a volume of a euclidean unit
ball). Similarly we can estimate other terms that appear in $\phi^* F$
and finally we can get

$$\mu(\bar{g},\frac{1}{2}) > (2\pi)^{-n}\int_M
e^{-\bar{f}_i}[|\nabla\bar{f}_i|^2 + R + \bar{f}_i - 2n] dV_{g_i} -
k\epsilon,$$ 
for big enough $i$. Therefore, $\mu(\bar{g},\frac{1}{2})
> \mu(g_i,\frac{1}{2}) - k\epsilon$ for $i \ge i_0$, for some constant $k$.
In a very similar manner as above we can get that $\mu(g_i,\frac{1}{2}) > 
\mu(g_i, \frac{1}{2}) - k\epsilon$ for $i \ge i_0$, since $f_i$ satisfy

$$2\Delta f_i - |\nabla f_i|^2 + R(g_i) + f_i - 2n = \mu(g_i,\frac{1}{2}),$$
where $C_1 \le \mu(g_i, \frac{1}{2}) \le C_2$.

\end{proof}

\begin{remark}
The previous proposition holds for a Riemannian manifold of an
arbitrary dimension.
\end{remark}

From now on restrict ourselves to a case of a K\"ahler-Ricci flow in
complex dimension $2$. Tian showed in \cite{tian2004} that in complex
dimension $2$ we do not need any curvature assumptions to show that if
$g(t)$ is a K\"ahler-Ricci flow then for every sequence $t_i\to\infty$
there exists a subsequence such that $(M,g(t_i+t))\to (Y,\bar{g}(t))$,
in the orbifold sense, where $(Y,\bar{g}(t))$ are the orbifolds with
only finitely many singularities. Our goal is to prove Theorem
\ref{theorem-theorem_soliton}, i.e.

\begin{theorem}
Let $(g_{k\bar{j}})_t = g_{j\bar{k}} - R_{j\bar{k}} =
\partial_j\bar{\partial}_k u$ be a K\"ahler-Ricci flow on a
$2$-dimensional, complex, K\"ahler manifold. Then for every sequence
$t_i\to\infty$, there exists a subsequence so that
$(M,g(t_i+t))\to(Y,\bar{g}(t))$ and $\bar{g}(t)$ is a K\"ahler-Ricci
soliton, i.e.

$$g_{i\bar{j}} - R_{i\bar{j}} = \partial_i\partial_{\bar{j}}\bar{u},$$
where $u_{ij} = u_{\bar{i}\bar{j}} = 0$ and $\bar{u}(t)$ is a
minimizer for $\mathcal{W}$ with respect to $\bar{g}(t)$.
\end{theorem}

\begin{proof}

Fix $r > 0$. Then the pullbacks of metrics $g(t_i + t)$ converge
uniformly to $\bar{g}(t)$ on $Y\backslash B_{\bar{g}(t)}(p,r/3)\times
[0,A]$, for any finite $A$. Therefore, there exists $i_0$ so that
$|\ric|(g(t_i+t)) \le C = C(r)$ and $|\rem|(g(t_i+t)) \le C$ for all
$i \ge i_0$ on $\phi_i(Y\backslash B_{\bar{g}(t)}(p,r/3))$, where
$\phi_i$ are diffeomorphisms from a definition of convergence and for
all $t\in [0,A]$.

(*)First of all we will make an appropriate choice of cut off
functions. Choose small $\delta < A$ so that

$$B_{g(t_i+s)}(p_i,r/2)\subset B_{g(t_i)}(p_i,r),$$
and
$$B_{g(t_i)}(p_i,2r)\subset B_{g(t_i+s)}(p_i,5r/2),$$ for all
$s\in[0,\delta]$, all $r > 0$ and all $i\ge i_0$. This is possible, 
since the uniform bounds on the Ricci curvatures give a good control
over the distances, the sizes of balls and the norms of vectors at
different times. This control does not depend on either the choice of
a point or time.

Let $\xi_i$ be a sequence of cut off functions such that $\xi_i = 1$
outside the ball $B_{g(t_i)}(p_i,2r)$, $\xi_i = 0$ in
$B_{g(t_i)}(p_i,r)$ and the norms of the derivatives of $\xi$ in
metrics $g(t_i)$ are uniformly bounded in $i$ by a constant that
dependes on $r$. We can decrease $\delta$ if necessary, so that for
$i\ge i_0$, norms of the derivatives of $\xi_i$ in metrics $g(t_i+s)$
are uniformly bounded in $i$ and $s\in [0,\delta]$. Then we have that
our cut off functions $xi_i$ acutually satisfy

\begin{displaymath}
\xi_i = \left\{\begin{array}{ll}
0 & \textrm{if $x\in B_{g(t_i+s)}(p_i,r/2)$}\\
1 & \textrm{outside $B_{g(t_i+s)}(p_i,5r/2)$}
\end{array} \right.
\end{displaymath}
for all $s\in [0,\delta]$ and all $i\ge i_0$.

For every $t$, choose a minimizer for Perelman's
functional $\mathcal{W}$ with respect to $g(t)$. Denote it by
$f_t$. Flow it backwards. Let $u_t = e^{-f_t}$. Let $t_i' = t_i +
\delta$. Then

$$\frac{d}{ds} u_{t_i'}(s) = -\Delta u_{t_i'}(s) + (n-R)u_{t_i'}(s).$$
Following Perelman's computation we get

$$\frac{d}{ds}\mathcal{W}(g(t_i+s),f_{t_i'}(s),1/2) = (2\pi)^{-n}\int_M
|R_{k\bar{j}} + \nabla_k\bar{\nabla}_j f_{t_i'} -
g_{k\bar{j}}|^2(g(t_i+s) dV_{g(t_i+s)}.$$
If we integrate it over $s\in[0,\delta]$, we get

\begin{eqnarray*}
\mu(g(t_i+\delta),1/2) - \mu(g(t_i),1/2) &\ge&
\mathcal{W}(g(t_i'),f_{t_i'},1/2) -
\mathcal{W}(g(t_i),f_{t_i'}(t_i),1/2)\\ &=&
(2\pi)^{-n}\int_0^{\delta}\int_M |\nabla\bar{\nabla}(f_{t_i'} -
u)|^2(t_i+s) dV_{g(t_i+s)}.
\end{eqnarray*}
$\mu(g(t),\frac{1}{2}) \le \mathcal{W}(g(t),u(t),1/2) =
(2\pi)^{-n}\int_M e^{-u}(|\nabla u|^2 + R + u - 2n) dV_{g(t)} \le C$,
because of Theorem \ref{theorem-theorem_useful}. Combining this with
Perelman's montonicity formula give that there exists a finite
$\lim_{t\to\infty}\mu(g(t),1/2)$ and therefore

$$\lim_{i\to\infty} |\nabla\bar{\nabla}(f_{t_i'} - u)|(t_i+s) = 0,$$
for almost all $x\in M$ and almost all $s\in [0,\delta]$.

Moreover,

$$\frac{d}{ds}\xi_i u_{t_i'}(s) = -\xi_i\Delta u_{t_i'}(s) +
\xi_i(n-R)u_{t_i'}(s),$$
for $s\in [0,\delta]$. As in \cite{natasa2004} we can get that

$$|u_{t_i'}(t_i+s)|_{C^{2,\alpha}(M\backslash B_{g(t_i+s)}(p_i,5r/2))}
\le C(r),$$ where $C(r)$ is a constant that depends on $r > 0$ which
we have fixed at the beginning. Since the curvature is uniformly
bounded on $M\backslash B_{g(t_i+s)}(p_i,5r/2)$ for all $i\ge i_0$ and
all $s\in[0,\delta]$, we have also that

$$|u(t_i+s)|_{C^{2,\alpha}(M\backslash B_{g(t_i+s)}(p_i,5r/2))} \le
C(r).$$ We can now extract a subsequence $t_i$ such that
$u_{t_i'}(t_i+s)\stackrel{C^{2,\alpha}}{\to}\bar{u}_1(s)$ and
$u(t_i+s)\stackrel{C^{2,\alpha}}{\to}\bar{u}(s)$ , on $Y\backslash
B_{\bar{g}(s)}(p,5r/2)$, uniformly in $s\in[0,\delta]$.  Because of
those $C^{2,\alpha}$ estimates, we have that

$$\nabla\bar{\nabla}(\bar{u} - \bar{f}) = 0,$$ i.e. $\Delta \bar{f} =
\Delta\bar{u}$ on $Y\backslash B_{\bar{g}(s)}(p,5r/2)$, where $\bar{u}
= e^{-\bar{f}}$.

We can consider instead of $r > 0$ a sequence $r_k = 1/k \to 0$.
Choosing a subsequence of a subsequence of a sequence $t_i$ for each
$k$-th step, using the uniqueness of a limit and a diagonalization
method we can extract a subsequence $t_i$ such that
$u_{t_i'}(t_i+s)\stackrel{C^{2,\alpha}}{\to}\bar{f}(s)$ and
$u(t_i+s)\stackrel{C^{2,\alpha}}{\to}\bar{u}(s)$ uniformly on compact
subsets of $Y\backslash\{p\}\times [0,\delta]$. We also know that
$\Delta\bar{f}(s) = \Delta\bar{u}(s)$ on $Y\backslash\{p\}$ and
therefore $\bar{f}(s) = \bar{u}(s)$, since both of them satisfy the
same integral normalization condition

$$\int_Y e^{-\bar{f}}dV_s = \int_Y e^{-\bar{u}}dV_s = (2\pi)^n,$$
and since $Y$ is a compact orbifold. 
Because of the uniform convergence on compact sets we have that
$\bar{u}(s) = \bar{f}(s)$ satisfy

$$\frac{d}{ds}\bar{f}(s) = -\Delta \bar{f} + |\nabla\bar{f}|^2 - R + n.$$
On the other hand $\bar{u}(s)$ satisfy

$$\frac{d}{ds}\bar{u} = \Delta\bar{u} + \bar{u} + a,$$
$$\Delta\bar{u} = n - R.$$
Therefore,

$$\frac{d}{ds}\bar{u}(s) = |\nabla\bar{u}|^2(s),$$
i.e.
$$\Delta\bar{u} - |\nabla\bar{u}|^2 + \bar{u} = -a.$$ $f_{t_i'}(t_i')
= f_{t_i'}$ is a minimizer for $\mathcal{W}$ with respect to $g(t_i')$
and it converges to a minimizer for $\mathcal{W}$ with respect to a
limit metric $\bar{g}(\delta)$, which is a consequence of Proposition
\ref{proposition-proposition_mu} (the arguments for this are similar
as in \cite{natasa2004}). Therefore, $\bar{u}(\delta)$ is a minimizer
for $\mathcal{W}$ with respect to $\bar{g}(t)$.

\begin{claim}
\label{claim-claim_a}
$\bar{a}(s) = a$, for all $s\in[0,\delta]$, where $\bar{a}(s) =
-(2\pi)^{-n}\int_Y \bar{u}e^{-\bar{u}}dV_s$.
\end{claim}

\begin{proof}

We can repeat everything that we have done before, replacing $\delta$
with any $t\in(0,\delta]$ to conclude that $\bar{u}(t)$ is a minimizer
for $\mathcal{W}$ with respect to $\bar{g}(t)$ for $t\in(0,\delta]$.

$$\mu(\bar{g}(t), 1/2) = (2\pi)^{-n}\int_Y (|\nabla\bar{u}|^2 -
\Delta\bar{u} + \bar{u} -n)e^{-\bar{u}} dV_{\bar{g}(s)} = -\bar{a}(s)
- n.$$ 
On the other hand there exists a finite $\lim_{t\to\infty}\mu(g(t),1/2)$
and therefore $\mu(\bar{g}(t),1/2) = \mu(\bar{g}(s),1/2)$ 
for all $s,t\in(0,\delta]$, i.e. $\bar{a}(t) = a$ for some 
constant $\bar{a}$ for all $s\in(0,\delta]$. Because of the continuity,
$\bar{a}(t) = a$ for all $t\in[0,\delta]$.

\end{proof}

\begin{claim}
\label{claim-claim_soliton}
$\bar{u}_{ij} = \bar{u}_{\bar{i}\bar{j}} = 0$ for all
$t\in[0,\delta]$.
\end{claim}

\begin{proof}

$$\Delta \bar{u} - |\nabla \bar{u}|^2 + \bar{u} = -a.$$ 
From here we see $(\Delta \bar{u})^2 = \Delta\bar{u} |\nabla\bar{u}|^2
- \bar{u}\Delta\bar{u} -a\Delta\bar{u}$, i.e.

\begin{equation}
\label{equation-equation_help1}
\int_M \Delta\bar{u}(|\nabla\bar{u}|^2 - \Delta\bar{u}) dV_{s} =
\int_M \bar{u}\Delta\bar{u} = -\int_M |\nabla\bar{u}|^2 dV_{s}
\end{equation}
We have the following evolution equations:

$$\frac{d}{dt} |\nabla\bar{u}|^2 = \Delta|\nabla\bar{u}|^2 -
|\nabla\nabla\bar{u}|^2 - |\nabla\bar{\nabla} \bar{u}|^2 +
|\nabla\bar{u}|^2.$$

$$\frac{d}{dt}\Delta\bar{u} = \Delta^2 \bar{u} + \Delta\bar{u} -
|\nabla\bar{\nabla}\bar{u}|^2.$$
If we subtract the second equation from the first one and if we
integrate what we get along $M$, we get

$$\int_M \frac{d}{dt}(|\nabla\bar{u}|^2 - \Delta\bar{u}) = -\int_M
|\nabla\nabla \bar{u}|^2 + \int_M |\nabla\bar{u}|^2,$$
i.e.
\begin{eqnarray}
\label{equation-equation_help2}
&&\frac{d}{dt}\int_M (|\nabla\bar{u}|^2 - \Delta\bar{u})dV_{\bar{g}} =\\ 
&=&-\int_M |\nabla\nabla \bar{u}|^2 + \int_M |\nabla\bar{u}|^2 + \int_M
(|\nabla\bar{u}|^2 - \Delta\bar{u})\Delta\bar{u} dV_{\bar{g}(t)} \nonumber\\
&=& -\int_M |\nabla\nabla \bar{u}|^2 dV_{h(t)} \nonumber,
\end{eqnarray}
where we have used the equation (\ref{equation-equation_help1}). 
On the other hand, we have that $\frac{d}{dt}\int_M 
(|\nabla\bar{u}|^2 - \Delta\bar{u})dV_{t} = \frac{d}{dt}
\int_M (\bar{u} + a)$.
Since $\frac{d}{dt}\bar{u}(t) = |\nabla\bar{u}|^2$, we have

$$\frac{d}{dt}\int_M\bar{u} dV_{\bar{g}(t)} = \int_M
(|\nabla\bar{u}|^2 + \bar{u}\Delta\bar{u}) dV_{\bar{g}(t)} = 0.$$ By
claim \ref{claim-claim_a}, $\bar{a}(t) = a$ is independent of $t$ and
volume is fixed along the flow. These imply that $\frac{d}{dt}\int_M
(\bar{u} + a) dV_{\bar{g}(t)} = 0$.
From equation (\ref{equation-equation_help2}) we get

$$\int_M |\nabla\nabla\bar{u}|^2 dV_{\bar{g}} = 0,$$
i.e. $\bar{u}_{ij} = \bar{u}_{\bar{i}\bar{j}} = 0$.

\end{proof}

Claim \ref{claim-claim_soliton} tells us that $\bar{u}(t)$ comes from
a holomorphic vector field on $Y$, i.e. $\bar{g}(t)$ is a
K\"ahler-Ricci soliton for $t\in[0,\delta]$. Since our $\delta$
depends on the uniform bouund on the Ricci curvature and not on an
initial time, we can apply the above proof of Theorem
\ref{theorem-theorem_soliton} to a sequence $t_i+\delta$ instead of a
sequence $t_i$ (with the same choice of $\delta$ as in $(*)$ above) to
conclude that there exists a subsequence $t_i$ such that $g(t_i+t)\to
\bar{g}(t)$ and $u(t_i+t)\stackrel{C^{2,\alpha}}{\to}\bar{u}(t)$
uniformly on compact subsets of $Y\backslash\{p\}\times[0,2\delta]$
and such that $(g_{k\bar{j}})_t = g_{k\bar{j}} - R_{k\bar{j}} =
\partial_k\partial_{\bar{j}}\bar{u}$ and $\bar{u}_{kj}(t) =
0$. Continuing this process and diagonalizing the sequence $t_i$ we
will get a subsequence $t_i$ so that
$(M,g(t_i+t))\to(Y,\bar{g}(t))$. The convergence is uniform on compact
subsets of $Y\backslash\{p\}\times [0,\infty)$ and $\bar{g}(t)$ is a
K\"ahler-Ricci soliton for all times $t$.

\end{proof}

\end{section}

\end{document}